\documentclass[journal, paper]{IEEEtran}

\usepackage{amsmath}
\usepackage{amsfonts}
\usepackage{graphics}
\usepackage{pstricks}
\usepackage{array}
\usepackage{float}
\usepackage{epsf}
\usepackage{graphicx}
\usepackage{multirow}
\usepackage{cite}
\usepackage{acronym}
\usepackage{comment}
\usepackage{upgreek}
\usepackage{subcaption}
\usepackage{nicefrac} 
\usepackage[font=small]{caption}

\newcommand{\bfg}[1]{\boldsymbol{#1}}
\newcommand{\bft}[1]{\hat{\boldsymbol{#1}}}
\newcommand{\bfm}[1]{{|\boldsymbol{#1}|}}
\newcommand{\bfp}[1]{\boldsymbol{#1}'}
\newcommand{\bfpp}[1]{\boldsymbol{#1}''}
\newcommand{\bfppp}[1]{\boldsymbol{#1}'''}
\newcommand{\bfd}[1]{\dot{\boldsymbol{#1}}}
\newcommand{\bfdd}[1]{\ddot{\boldsymbol{#1}}}
\newcommand{\bfddd}[1]{\dddot{\boldsymbol{#1}}}
\newcommand{\bfb}[1]{\boldsymbol{\rm #1}}
\newcommand{\flux}{\bfg \varphi}
\newcommand{\ii}{\imath}

\newcommand{\fv}[1]{\boldsymbol{\rm f}_{#1}}
\newcommand{\e}[1]{\boldsymbol{\rm e}_{\rm #1}}
\newcommand{\ei}[1]{\boldsymbol{\rm e}_{#1}}
\newcommand{\eh}[1]{\tilde{\boldsymbol{\rm e}}_{#1}}
\newcommand{\sw}{\bfg \omega}

\newcommand{\T}{{\scriptstyle \rm \bfb T}}
\newcommand{\B}{{\scriptstyle \rm \bfb B}}
\newcommand{\N}{{\scriptstyle \rm \bfb N}}
\newcommand{\Td}{\dot{\scriptstyle \rm \bfb T}}
\newcommand{\Bd}{\dot{\scriptstyle \rm \bfb B}}
\newcommand{\Nd}{\dot{\scriptstyle \rm \bfb N}}
\renewcommand{\th}[1]{\theta_{#1}}
\newcommand{\ta}{\theta_a}
\newcommand{\wa}{\omega_a}
\newcommand{\wo}{\omega_o}
\renewcommand{\wr}{\omega_{r}}
\newcommand{\ws}{\omega_{s}}
\newcommand{\wk}{\omega_{\kappa}}
\newcommand{\wt}{\omega_{\tau}}
\renewcommand{\wp}{\omega_{\scriptscriptstyle \rm P}}
\newcommand{\tp}{\theta_{\scriptscriptstyle \rm P}}

\newcommand{\w}[1]{\omega_{#1}}
\renewcommand{\P}{\bfb P}
\newcommand{\F}{\bfb F}
\newcommand{\A}{\bfb A}
\newcommand{\C}{\bfb C}
\newcommand{\Tr}[1]{\hat{\bfb #1}}
\newcommand{\W}[1]{\bfb \Omega_{\scriptscriptstyle \rm #1}}
\newcommand{\Wt}[1]{\hat{\bfb \Omega}_{\scriptscriptstyle \rm #1}}
\newcommand{\Inv}[1]{\bar{\bfb #1}}
\newcommand{\vabc}[1]{\bfg #1_{abc}}
\newcommand{\vdqo}[1]{\bfg #1_{dqo}}
\newcommand{\vtnb}[1]{\bfg #1_{\scriptscriptstyle \rm TNB}}
\newcommand{\proj}[2]{\varpi_{#1}(#2)}

\begin{document}

\title{The Frenet Frame as a Generalization\\of the Park Transform}

\author{Federico~Milano,~\IEEEmembership{IEEE Fellow}\vspace{-4mm}%
  \thanks{F.~Milano is with School of Electrical and Electronic
    Eng., University College Dublin, Belfield, D04V1W8, Ireland.
    E-mail: federico.milano@ucd.ie}%
  \thanks{This work is supported by the Sustainable Energy Authority
    of Ireland (SEAI) under project FRESLIPS, Grant No.~RDD/00681; and
    by the European Commission under project EdgeFLEX, Grant
    No.~883710.}%
}

\maketitle

\begin{abstract}
  The paper proposes a generalization of the Park transform based on
  the Frenet frame, which is a special set of coordinates defined in
  differential geometry for space curves.  The proposed geometric
  transform is first discussed for three dimensions, which correspond
  to the common three-phase circuits.  Then, the expression of the
  time derivative of the proposed transform is discussed and the
  Frenet-Serret formulas and the Darboux vector are introduced.  The
  change of reference frame and its differentiation based on Cartan's
  moving frames and attitude matrices are also described.  Finally,
  the extension to circuits with more than three phases is presented.
  The features of the Frenet frame are illustrated through a variety
  of examples, including a case study based on the IEEE 39-bus system.
\end{abstract}


\IEEEpeerreviewmaketitle

\vspace{-2mm}

\section{Notation}

Scalars are indicated with Italic font, e.g.~$x$, whereas vectors and
matrices are indicated in bold face, e.g.~$\bfg x = (x_1, x_2, x_3)$.
Vectors have order 3, unless otherwise indicated.

\vspace{2mm}

\subsubsection*{Scalars}

\begin{itemize}[\IEEEiedlabeljustifyl \IEEEsetlabelwidth{Z} \labelsep 0.9cm]
\item[$s$] length of a curve
\item[$t$] time
\item[$v$] voltage magnitude
\item[$\alpha$] $2\pi/3$ rad
\item[$\beta$] $2\pi/6$ rad
\item[$\delta$] synchronous machine rotor angle
\item[$\theta$] voltage phase angle
\item[$\kappa$] curvature
\item[$\varpi$] projection operator
\item[$\tau$] torsion
\item[$\chi$] generalized curvature
\item[$\omega$] angular frequency
\end{itemize}

\vspace{2mm}

\subsubsection*{Vectors}

\begin{itemize}[\IEEEiedlabeljustifyl \IEEEsetlabelwidth{Z} \labelsep 0.9cm]
\item[$\bfg 0$] null vector %
\item[$\B$] binormal vector of the Frenet frame %
\item[$\ei{i}$] $i$-th vector of an orthonormal basis
\item[$\fv{i}$] $i$-th vector of a generalized Frenet frame
\item[$\bfg \imath$] current vector
\item[$\N$] normal vector of the Frenet frame %
\item[$\bfb r$] Darboux angular momentum vector %
\item[$\T$] tangent vector of the Frenet frame %
\item[$\bfg v$] voltage vector 
\item[$\bfg \phi$] magnetic flux vector
\end{itemize}

\vspace{2mm}

\subsubsection*{Matrices}

\begin{itemize}[\IEEEiedlabeljustifyl \IEEEsetlabelwidth{Z} \labelsep 0.9cm]
\item[$\A$] attitude matrix of Cartan's moving frame %
\item[$\C$] cylindrical frame %
\item[$\F$] matrix of the Frenet frame %
\item[$\P$] matrix of the Park transform %
\item[$\bfb \Psi$] generalized matrix to change coordinates %
\item[$\bfb \Omega$] rotation matrix %
\end{itemize}

\vspace{2mm}

\subsubsection*{Vector and Matrix Operations}

\begin{itemize}[\IEEEiedlabeljustifyl \IEEEsetlabelwidth{Z} \labelsep 1.8cm]
\item[$a, \bfm a$] vector magnitude
\item[$a', \bfg a', \A'$] derivative of a scalar/vector/matrix w.r.t. $t$
\item[$\dot{a}, \dot{\bfg a}, \dot \A$] derivative of a
  scalar/vector/matrix w.r.t. $s$
\item[$\bft{a}, \Tr{A}$] transpose of a vector/matrix
\item[$\Inv{A}$] inverse of a matrix
\end{itemize}

\vspace{-2mm}

\section{Introduction}
\label{sec:intro}

\subsection{Motivation}

The Park transform is the most important transform utilized in power
system transient stability analysis and control.  Originally
formulated by Park for the two-reaction theory of synchronous machines
\cite{5055275}, this transform has found applications in the control
of induction machines and, more recently, of converter-interfaced
devices.  In simulations, the Park transform is also a fundamental
tool for the implementation of the devices that compose the grid
\cite{FDF:2020} and for the modeling on power electronic converters
\cite{Iravani:2010}.  Despite its relevance, there are not many
attempts to generalize the Park transform, except for some extensions
to multi-phase circuits \cite{GABER1988203, 4454446}, nor to overcome
its intrinsic idiosyncrasies or better understand its geometric
properties.  This work addresses precisely these issues through the
theory provided by differential geometry.

\vspace{-2mm}
\subsection{Literature Review}

Since the introduction of phasors by Steinmetz at the end of the 19th
century \cite{Steinmetz:1897}, domain and coordinate transformations
are a common practice for the analysis of electrical circuits,
electrical machines and for power system modeling, analysis and
control.  Apart from the aforementioned Park transform, well-known
transforms are the harmonic analysis through Fourier series
\cite{Das:2015}, Fortescue symmetrical component theory
\cite{Fortescue:1918}, the forward-backward transform \cite{Ku:1929},
the Clarke transform \cite{Clarke:1943}, and more recently, dynamic
phasor analysis \cite{871734, 6687278}.  Except for the Fourier
analysis, the transforms above can be represented as 3-by-3 matrices
when applied to three-phase circuits.  Begin homomorphisms, it is also
possible to find the conversion matrices from one transform to another
(see, e.g., \cite{Hancock:1974}).  In the same vein, this work aims at
discussing a generalization of the concept of change of coordinates
and considers the most general case, i.e., the case for which the axis
of the coordinates are time dependent.

The proposed generalization is inspired by the literature on
differential geometry and, in particular, by the Frenet frame, the
Frenet-Serret formulas and the theory of moving frames developed by
Cartan \cite{Stoker, ONeill}.  These have found several applications
in mechanics, e.g., just to cite a recent relevant one, in the area of
autonomous vehicle driving \cite{2003lapierre, 2010werling}.  On the
other hand, the Frenet frame has been only very recently considered
for circuit and power systems analysis \cite{freqfrenet}.  This was
given raise by the geometrical interpretation of the frequency
developed by the author in \cite{freqgeom}.

References \cite{freqfrenet} and \cite{freqgeom} assume that the
instantaneous values of electrical quantities such as voltages and
currents of three- or multi-phase circuits are \textit{vectors} in a
given coordinate system.  This idea was already developed in the past
in the context of the instantaneous power theory \cite{4450060503,
  1308315, 5316097, IPT}.  These works define the active and reactive
power as the dot and cross (or wedge in multi-phase circuits)
products, respectively, of voltage and currents.  A variety of recent
works with same starting point have focused on the analysis of
electric circuits using the formalism provided by geometric algebra
\cite{math9111295, 2016barry, 2018ishihara, 2015talebi}.

The main difference of \cite{freqfrenet} and \cite{freqgeom} with
respect to the more conventional theory on instantaneous power is the
hypothesis that the voltage (current) vector is the \textit{speed} of
a curve, represented by the magnetic flux (electric charge).  This
interpretation moves the focus from ``algebra'' to ``calculus'' and
allows interpreting the dynamic behavior of the voltage and current in
terms of the ``invariants'' of differential geometry, such as curve
length, curvature and torsion.

\subsection{Contributions}

This work elaborates on the results of \cite{freqfrenet} and
determines similarities and differences between the Park transform and
the Frenet frame.  The novel contributions are the following.
\begin{itemize}
\item The derivation of the formal conditions under which the Park
  transform and its time derivative is equivalent to the Frenet frame
  and to the Frenet-Serret formulas.
\item An application of Cartan's moving frames that allows interfacing
  the local Frenet frame of each device connected to the grid to the
  reference frame of the grid itself.
\item An application of generalized $n$-dimensional Frenet frame to
  multi-phase systems, i.e., systems with more than three phases.
\item A thorough example-based discussion of the added value of the
  Frenet frame compared to the Park transform for the study of the
  dynamic performance of power systems.
\end{itemize}

\subsection{Organization}

The remainder of the paper is organized as follows.  Section
\ref{sec:park} recalls the definition of the Park transform and its
time derivative.  Section \ref{sec:frenet} introduces geometric
calculus, provides the definitions of curve length, curvature and
torsion, recalls the Frenet frame, the Frenet-Serret formulas and
Cartan's moving frames and provides the formulas of the generalized
curvatures in $n$ dimensions.  Section \ref{sec:geometry} combines the
definitions provided in the previous sections and defines the
conditions under which the Park transform is a special case of the
Frenet frame.  The interconnection of local Frenet frames of
electrical devices with the reference frame of the grid is also
discussed in Section \ref{sec:geometry} by means of Cartan's moving
frames.  Section \ref{sec:examples} illustrates the differences
between the Park transform and the Frenet frame through a series of
examples in three and six dimensions as well as the IEEE 39-bus
system.  Section \ref{sec:conc} draws conclusions and outlines future
work.

\section{Park Transform}
\label{sec:park}

The Park transform projects the phase components $abc$ of a
three-phase electrical quantity onto a $dqo$ frame, where the axes $d$
and $q$ rotate at angular speed $\wp$.  In his original formulation,
Park aimed at preserving the magnitudes of the transformed quantities
and did not define a power invariant transformation.  For the
developments presented in this paper, however, it is more convenient
to retain power invariance.

The formulation of the $dqo$-transform utilized in this paper is as
follows:
\begin{equation}
  \label{eq:abc2dqo}
  \bfg v_{dqo}(t)  = \P(\wp) \, \bfg v_{abc}(t) \, ,
\end{equation}
where $\bft v_{dqo} = [v_d, v_q, v_o]$,
$\bft v_{abc} = [v_a, v_b, v_c]$, and
\begin{equation}
  \label{eq:park}
    \P(\wp) = 
    \sqrt{\frac{2}{3}}
    \begin{bmatrix}
      \sin(\tp) & \sin(\tp - \alpha) & \sin(\tp + \alpha) \\
      \cos(\tp) & \cos(\tp - \alpha) & \cos(\tp + \alpha) \\
      \nicefrac{1}{\sqrt{2}} & \nicefrac{1}{\sqrt{2}} & \nicefrac{1}{\sqrt{2}}
    \end{bmatrix}
    \normalsize ,  
\end{equation}
where:
\begin{equation}
  \label{eq:tp}
  \begin{aligned}
    \tp(t) &= \int_0^t \wp(r) \, dr  +
    \theta_{{\scriptscriptstyle \rm P},0} \, .
  \end{aligned}
\end{equation}
Note that $\wp$ in \eqref{eq:tp} does not have to be constant.  The
power invariance of the matrix $\P$ in \eqref{eq:park} refers to the
fact that if $\bfg v_{abc}$ and $\bfg \ii_{abc}$ are the voltage and
current at a given point of a three-phase circuits, then the
instantaneous power is unchanged for the same voltage and current
transformed in $dqo$ coordinates:
\begin{equation}
  \bft v_{abc}(t) \, \bfg \ii_{abc}(t) =
  \bft v_{dqo}(t) \, \bfg \ii_{dqo}(t) \, ,
\end{equation}
%
%
%
This property descends from the fact that $\P$ is orthonormal,
i.e., its transpose is equal to its inverse:
\begin{equation}
  \label{eq:ti}
  \Tr{P}(\wp) = \Inv{P}(\wp) \, .
\end{equation}
Equation \eqref{eq:ti} can be readily proved, as follows:
%
\begin{equation*}
  \begin{aligned}
    \bft v_{dqo} \bfg \ii_{dqo} 
    &= \widehat{(\P \, \bfg v_{abc})} \P \, \bfg \ii_{abc} \\
    &= \bft v_{abc} \, \Tr{P} \P \, \bfg \ii_{abc}  \\
    &= \bft v_{abc} \, \Inv{P} \P \, \bfg \ii_{abc} 
    = \bft v_{abc} \bfg \ii_{abc} \, ,
  \end{aligned}
\end{equation*}
where the dependencies on time and on $\wp$ have been dropped for
simplicity.

The choice for the angular speed $\wp$ of the Park transform depends
on the device and the application.  For synchronous machines, $\wp$ is
generally chosen as the rotor angular speed of the machine itself,
namely $\wp = \omega_r = \delta_r'$.  This allows simplifying the
equations of the machine and rewriting rotor quantities and equations
as they were a dc circuit.  This is also the motivation for the
original two-reaction theory developed by Park.  For all other
devices, however, including induction machines, and in general for
transient stability analysis studies of interconnected systems, it is
chosen $\wp = \wo$, namely, the constant synchronous reference angular
frequency of the grid, e.g., $\wo = 2 \pi \, 60$ rad/s.

\subsection{Time Derivative of $dqo$-Axis Voltages}
\label{sub:dqodot}

The time derivative of a $dqo$-axis voltage is given by:
\begin{equation}
  \label{eq:vdqodot}
  \begin{aligned}
    \bfg v'_{dqo} &= \big (\P \, \bfg v_{abc} \big )' \\
    &= \P' \bfg v_{abc} + \P \, \bfg v'_{abc} \\
    &= \P' \, \Inv{P} \bfg v_{dqo} + \P \, \bfg v'_{abc} \, ,
  \end{aligned}
\end{equation}
where the dependency of $\P$ on $\wp$ has been dropped for simplicity.
Let us define:
\begin{equation}
  \label{eq:Pw}
  \W{P} = \P' \, \Inv{P} = \P' \, \Tr{P} =
  \begin{bmatrix}
    0 & \wp & 0 \\
    -\wp & 0 & 0 \\
    0 & 0 & 0
  \end{bmatrix} ,
\end{equation}
which is a skew-symmetric matrix, i.e., $\Wt{P} = -\W{P}$.  Then,
\eqref{eq:vdqodot} can be rewritten as:
\begin{equation}
  \label{eq:vabcdot}
  \begin{aligned}
    \P \, \bfg v'_{abc}
    &= \bfg v'_{dqo} - \W{P} \, \bfg v_{dqo} \\
    &= \bfg v'_{dqo} + \Wt{P} \, \bfg v_{dqo} \, .
  \end{aligned}
\end{equation}
Equation \eqref{eq:vabcdot} shows that the Park transform of the
derivative of $\bfg v_{abc}$ consists of two terms.  The first term is
the time derivative of the Park-transformed voltage $\bfg v_{dqo}$,
which represents a translation.  The second term is given by product
of $\W{P}$ and $\bfg v_{dqo}$, which is due to the rotation of the
Park $dq$-axis.  The second term is null only if $\wp = 0$, which is
the choice that leads to the Clarke transform.


\subsection{Interface between $dqo$-Axes Rotating at Different Speeds}
\label{sub:C}

We consider the relevant case of the interface of the stator
terminal-bus synchronous generator with Park transform $\P(\wr)$ with
the grid, the Park transform of which is $\P(\ws)$.  Since the
zero-axis of the Park transform does not rotate with $\wp$, the
interface between machine and network consists in a rotation in the
$dq$-plane, as follows:
\begin{equation}
  \begin{aligned}
  \label{eq:G2N}
  \vdqo{v}^N (t)
  &= \P(\ws) \, \vabc{v}(t) \\
  &= \P(\ws) \, \Tr{P}(\wr) \, \vdqo{v}^G (t) \\
  &= \C \, \vdqo{v}^G (t) \, ,
\end{aligned}
\end{equation}
where $N$ and $G$ indicate ``network'' and ``generator'',
respectively, and matrix $\C$ is orthonormal and represents a change
of coordinates in a cylindrical frame.  In fact, developing the matrix
multiplication between $\P(\ws)$ and $\Tr{P}(\wr)$, one obtains:
\begin{equation}
  \label{eq:C}
  \C = 
  \begin{bmatrix}
    \cos(\theta_s -\delta_r) & -\sin(\theta_s -\delta_r) & 0 \\
    \sin(\theta_s -\delta_r) & \cos(\theta_s -\delta_r) & 0 \\
    0 & 0 & 1 \\
  \end{bmatrix} ,
\end{equation}
where $\delta_r = \int_0^t \wr(r) - \wo \, dr + \delta_{r,0}$ is the
rotor angular position of the machine and
$\theta_s= \int_0^t \ws(r) - \wo \, dr + \theta_{s,0}$ is the phase
angle of the terminal-bus voltage of the machine.  In practical
implementations, one can only know a relative value of the phase
angles.  For this reason $\delta_r$ and $\theta_s$ must be referred to
the same reference frame, which is generally chosen as rotating at the
constant reference angular speed $\wo$.  Another common choice for the
reference is the angular frequency of the center of inertia of the
system \cite{FDF:2020}.  This is utilized to avoid the drift of the
machine angles and bus voltage phase angles during the transients
following a large perturbation \cite{Fabozzi:2011}.


Finally, the time derivative at the interface between a generator and
the network can be obtained in a similar way as described in the
previous section, as follows:
\begin{equation}
  \label{eq:vdqodotGN}
  \begin{aligned}
    (\vdqo{v}^N)' &= \big (\P(\ws) \, \vabc{v} \big )' \\
    &= \P'(\ws) \vabc{v} + \P(\ws) \, \vabc{v}' \\
    &= \P'(\ws) \, \Tr{P}(\wr) \vdqo{v}^G + \P(\ws) \, \vabc{v}' \, ,
  \end{aligned}
\end{equation}
and substituting
\begin{equation}
  \vabc{v}' = \Tr{P}(\wr) (\vdqo{v}^G)' - \Tr{P}(\wr) \W{P}(\wr) \vdqo{v}^G \, ,
\end{equation}
one obtains:
\begin{equation}
  \label{eq:vdqodotGN2}
  \begin{aligned}
    (\vdqo{v}^N)' &= \C \, (\vdqo{v}^G)'  \\
    &+
    [\P'(\ws) \, \Tr{P}(\wr) - \C \, \W{P}(\wr)] \vdqo{v}^G \, .
  \end{aligned}
\end{equation}
Observing that:
\begin{equation}
  \begin{aligned}
    \P'(\ws) \, \Tr{P}(\wr)
    &= \P'(\ws) \, \Tr{P}(\ws) \, \P(\ws) \, \Tr{P}(\wr) \\
    &= \W{P}(\ws) \, \C \, ,
  \end{aligned}
\end{equation}
and that $\W{P}(\ws) \, \C = \C \, \W{P}(\ws)$, one
obtains:
\begin{equation}
  \label{eq:vdqodotGN3}
  \begin{aligned}
    (\vdqo{v}^N)' &= \C \, (\vdqo{v}^G)' \\
    &+
    \C \, [\W{P}(\ws) - \W{P}(\wr)] \vdqo{v}^G \, ,
  \end{aligned}
\end{equation}
where it is relevant to observe that the term
$\W{P}(\ws) - \W{P}(\wr)$ can be also obtained as:
\begin{equation}
  \label{eq:Cw}
  \C' \, \Tr{C} =
  \begin{bmatrix}
    0 & \ws - \wr & 0 \\
    \wr - \ws & 0 & 0 \\
    0 & 0 & 0
  \end{bmatrix} .
\end{equation}
The latter expression does not appear out of the blue.  It is the
consequence of a more general theory, i.e., Cartan's moving frames,
which is described in the next section.

\section{Frenet Frame of Space Curves}
\label{sec:frenet}

This section introduces the classical Frenet frame and the
Frenet-Serret formulas of space curves.  With this aim, it is relevant
to provide first some definitions.

The starting point is a curve in a three-dimensional space, say
$\bfg x = (x_1, x_2, x_3)$ or, equivalently:
\begin{equation}
  \bfg x = x_1 \, \e{1} + x_2 \, \e{2} + x_3 \, \e{3} \, ,
\end{equation}
where $(\e{1}, \e{2}, \e{3})$ is an orthonormal basis.
For the development given below, it is relevant to define two types of
products that can be done with three-dimensional vectors, namely the
dot product and and the cross product.  The \textit{dot product} of
two vectors returns a scalar, as follows:
\begin{equation}
  \label{eq:inner2}
  \bfg x \cdot \bfg y = \bft x \, \bfg y = x_1y_1 + x_2y_2 + x_3y_3 \, .
\end{equation}
The \textit{cross product} of two vectors returns a vector that is
orthogonal to the original vectors, as follows:
\begin{equation}
  \begin{aligned}
    \bfg x \times \bfg y
    &= 
      \left | 
      \begin{matrix}
        \e{1} & \e{2} & \e{3} \\
        x_1 & x_2 & x_3 \\
        y_1 & y_2 & y_3 \\
      \end{matrix} 
    \right | \, .
  \end{aligned}
\end{equation}
It is relevant to note that $\bfm x = \sqrt{\bfg x \cdot \bfg x}$ is
the magnitude of the vector, and $\bfg x \times \bfg x = \bfg 0$.

The length $s$ of the curve is defined as:
\begin{equation}
  s = \int_0^t \sqrt{(\bfp x \cdot \bfp x)} \, dt \, ,
\end{equation}
or, equivalently:
\begin{equation}
  \label{eq:s}
  s' = \sqrt{\bfp x \cdot \bfp x} = |\bfp x| \, .
\end{equation}
where
\begin{equation}
  \begin{aligned}
    \bfg x'
    =& \; (x_1 \, \e{1})' +
      (x_2 \, \e{2})' +
      (x_3 \, \e{3})' \\
    =& \; (x'_1 \, \e{1} +
      x'_2 \, \e{2} +
      x'_3 \, \e{3}) + \\
    & \; 
      (x_1 \, \e{1}' +
      x_2 \, \e{2}' +
      x_3 \, \e{3}') 
  \end{aligned}
\end{equation}
is the speed of the trajectory described by $\bfg x$.  For fixed
reference frames the terms $\e{i}'$ are null.  However, in this work,
it is of interest to discuss moving frames, i.e., sets of coordinates
for which the position of the axes of the coordinates vary in time.
Both the Park transform and the Frenet frame, which are introduced
below, are special cases of moving frames.

The length $s$ is a geometric invariant, i.e., its value does not
depend on the choice of the coordinates.  Neglecting relativistic
effects, also its time derivative, $s'$, is an invariant and has a
special role in differential geometry.  In particular, it is relevant
to define The derivative of $\bfg x$ with respect to $s$, which, using
the chain rule, can be written as:
\begin{equation}
  \label{eq:xdot}
  \bfd x = \frac{d \bfg x}{d s} =
  \frac{d\bfg x}{dt} \frac{dt}{ds} = \frac{\bfp x}{s'} =
  \frac{\bfp x}{|\bfp x|} \, .
\end{equation}
The vector $\bfd x$ has magnitude 1 and is tangent to the curve
$\bfg x$.  In the remainder of this paper, $\bfd x$, $\bfdd x$,
$\bfddd x$, etc.~indicate the derivatives of a vector with respect to
$s$, whereas $\bfg x'$, $\bfg x''$, $\bfg x'''$ indicate time
derivatives.

We have mentioned that invariants play a relevant role in differential
geometry as they are quantities independent form the choice of the
coordinates.  Yet, among all possible set of coordinates, there is
one, the Frenet frame, that has special properties.  This frame is
defined by the following three vectors:
\begin{equation}
\label{eq:TNB}
  \begin{aligned}
    &\T = \bfd x \, , \qquad
    &\N = \frac{\bfdd x}{|\bfdd x|} \, , \qquad
    &\B = \T \times \N \, ,
  \end{aligned}
\end{equation}
where $\T$, $\N$ and $\B$ are called tangent, normal and binormal
vectors, respectively.  The vectors in \eqref{eq:TNB} are orthonormal,
i.e.~$\T = \N \times \B$ and $\N = \B \times \T$, and satisfy the
following set of differential equations \cite{Stoker}:
\begin{equation}
  \label{eq:frenet3}
  \begin{aligned}
    \Td &= \kappa \, \N \, , \\
    \Nd &= -\kappa \, \T + \tau \, \B \, , \\
    \Bd &= - \tau \, \N \, ,
  \end{aligned}
\end{equation}
where $\kappa$ and $\tau$ are the \textit{curvature} and the
\textit{torsion}, respectively, which are given by:
\begin{equation}
  \label{eq:kappa}
  \kappa = |\bfdd x| = \frac{|\bfp x \times \bfpp x|}{|\bfp x|^3} \, , 
\end{equation}
and
\begin{equation}
  \label{eq:tau}
  \tau = \frac{\bfd x \cdot \bfdd x \times \bfddd x}{\kappa^2} =
  \frac{\bfp x \cdot \bfpp x \times \bfppp x}{|\bfp x \times \bfpp x|^2} \, .
\end{equation}
Both $\kappa$ and $\tau$ are geometric invariants, as the length $s$.

Equations \eqref{eq:frenet3} are known as Frenet-Serret equations and
have a key role in this paper.  Observe that \eqref{eq:frenet3} can be
rewritten as:
\begin{equation}
  \label{eq:frenet4}
  \dot {\F} =
  \begin{bmatrix}
    0 & \kappa & 0 \\
    -\kappa & 0 & \tau \\
    0 & -\tau & 0
  \end{bmatrix} \F \, ,
\end{equation}
where
\begin{equation}
  \label{eq:frenet5}
  \F =
  \begin{bmatrix}
    \hat{\T} \\ \hat{\N} \\ \hat{\B}
  \end{bmatrix} =
  \begin{bmatrix}
    {\scriptstyle T}_{1} & {\scriptstyle T}_{2} & {\scriptstyle T}_{3} \\
    {\scriptstyle N}_{1} & {\scriptstyle N}_{2} & {\scriptstyle N}_{3} \\
    {\scriptstyle B}_{1} & {\scriptstyle B}_{2} & {\scriptstyle B}_{3}    
  \end{bmatrix} .
\end{equation}
Recalling \eqref{eq:xdot}, the Frenet-Serret equations
\eqref{eq:frenet3} and, thus, \eqref{eq:frenet4}, can be rewritten,
using the chain rule, as:
\begin{equation}
  \label{eq:Fdot}
  \F' = \dot \F \, s' =
  \begin{bmatrix}
    0 & \wk & 0 \\
    -\wk & 0 & \wt \\
    0 & - \wt & 0
  \end{bmatrix} \F = \W{F} \, \F \, ,
\end{equation}
where $\wk = s' \kappa$ and $\wt = s' \tau$ have the dimension of an
angular frequency and are defined in \cite{freqfrenet} as
\textit{azimuthal frequency} and \textit{torsional frequency},
respectively.  Similarly to $\kappa$ and $\tau$, also $\wk$ and $\wt$
are geometric invariants, since they are products of invariants.

From \eqref{eq:Fdot} and the orthonormality of $\F$ descend that
\cite{ONeill}:
\begin{equation}
  \label{eq:Fw}
  \W{F} = \F' \, \Inv{F} = \F' \, \Tr{F} \, ,
\end{equation}
whose resemblance with \eqref{eq:Pw} is evident.


\subsection{Time Derivative and Darboux Vector}
\label{sub:darboux}

This section is dual to Section \ref{sub:C}, i.e., discusses the time
derivative of vectors transformed using the Frenet frame.  Let us define:
\begin{equation}
  \label{eq:vtnb}
  \vtnb{x} = \F \, \bfg x \, .
\end{equation}
Then, the time derivative of $\vtnb{x}$ is:
\begin{equation}
  \label{eq:vtnbdot1}
  \begin{aligned}
    \vtnb{x}' &= \F \, \bfp x + \F' \, \bfg x \\
    &= \F \, \bfp x + \F' \, \Inv{F} \, \vtnb{x}  \\
    &= \F \, \bfp x + \W{F} \, \vtnb{x}  \, ,
  \end{aligned}
\end{equation}
And, finally:
\begin{equation}
  \label{eq:vtnbdot}
  \F \, \bfp x = \vtnb{x}' + \Wt{F} \vtnb{x} \, ,
\end{equation}
which has the same structure as \eqref{eq:vabcdot}.  It is relevant to
note that:
\begin{equation}
  \bfb r  = \wt \, \T + \wk \, \B = \star ( \W{F} ) \, ,
\end{equation}
where $\bfb r$ is called \textit{Darboux vector} (or \textit{angular
  momentum vector}) and $\star$ is the Hodge star operator, i.e., an
isomorphism between vectors and matrices (bivectors).  Hence, one can
rewrite \eqref{eq:vtnbdot} as \cite{Stoker}:
\begin{equation}
  \label{eq:K}
  \F \, \bfp x = \vtnb{x}' - \bfb r \times \vtnb{x} \, ,
\end{equation}
where $\F \bfp x$ is the transform of the speed of the curve $\bfg x$
and can be interpreted as the time derivative of the transformed curve
$\vtnb{x}$ on the rotating frame defined by $\F$ plus a term that
depends on the rotation of the Frenet frame itself.  Differently from
the Park transform, however, the Frenet frame has two rotations: in
the plane $(\T, \N)$ with angular frequency $\wk$ and in the plane
$(\N, \B)$ with angular frequency $\wt$.

\subsection{Cartan's Moving Frames}

While the Darboux vector is useful to understand the geometrical
meaning of the rotation matrix $\W{F}$, it cannot be easily generalized
to an arbitrary set of orthonormal basis.  This generalization is due
to Cartan,\footnote{Cartan's notation utilizes the differential 1-form
  $\bfb d$ rather than the time derivative.  However, for sake of
  simplicity and consistency with the conventional vector-based
  notation utilized in power systems, differential forms are not used
  in this work.  The interested reader can find an introduction to
  Cartan's forms and moving frames in \cite{Needham:2021}.}  who
obtained the following general expression for an orthonormal matrix
$\A$:
\begin{equation}
  \label{eq:Aw0}
  \W{A} = \A' \, \Inv{A} = \A' \, \Tr{A} \, ,
\end{equation}
where $\W{A}$ has the following general structure:
\begin{equation}
  \label{eq:Aw1}
  \W{A} =
  \begin{bmatrix}
    0 & \omega_{12} & \omega_{31} \\
    -\omega_{12} & 0 &  \omega_{23} \\
    -\omega_{31} & -\omega_{23} & 0
  \end{bmatrix} ,
\end{equation}
which is a skew-symmetric matrix, namely $\Wt{A} = -\W{A}$.

In \cite{Needham:2021}, $\A$ is called \textit{attitude} (or
\textit{orientation}) {matrix}.  Thus, $\F$ is an attitude matrix for
which $\omega_{31} = 0$.  It is possible to demonstrate that the
information contained in $\W{F}$ is \textit{all} the information on
the rotation of the frame \cite{ONeill}.  Hence, $\F$ is not just any
attitude matrix, but the matrix that represents a particular (the
only, in fact) moving frame following the trajectory of the curve
$\bfg x$ for which $\omega_{31} = 0$ for all $t$.

Based on the definition above, $\P(\wp)$ is also an attitude matrix
for which $\omega_{12} = \wp$ and $\omega_{23} = \omega_{31} = 0$.
However, since for $\P$, $\omega_{23}$ is \textit{always} null, the
information provided by $\P$ is incomplete in general.  This point is
a key contribution of the paper and is further elaborated in Section
\ref{sec:geometry}.

\subsection{Extension to $n$-Dimensional Curves}
\label{sub:ndim}

So far, we have discussed three-dimensional vectors and three-phase
circuits.  However, the Frenet frame and Frenet-Serret formulas can be
extended to $n$ dimensions, and hence, to circuits with an arbitrary
number of phases.  The starting point are a set of independent
vectors, e.g., $(\bfd x, \bfdd x, \dots , \bfg x^{(n)})$.  Then, the
orthonormal basis can be constructed through the Gram-Schmidt process.
With this aim, let define first the \textit{projection operator} as:
\begin{equation}
  \proj{\bfg u}{\bfg w}
  = \frac{\bfg u \cdot \bfg w}{\bfg u \cdot \bfg u} \, .
\end{equation}
Then the unnormalized orthogonal vectors
$\eh{1}, \eh{2}, \dots, \eh{n}$ are obtained as:
\begin{equation}
  \begin{aligned}
    \eh{1} &= \bfd x \, , \\
    \eh{h} &= \bfg x^{(h)} - \sum_k^{h-1} \proj{\eh{k}}{\bfg{x}^{(h)}} \, , \quad h = 2, \dots, n-1 \, , 
  \end{aligned}
\end{equation}
Then, the normalized orthonormal $n-1$ vectors are given by:
\begin{equation}
  \begin{aligned}
    \fv{i} = \frac{\eh{i}}{|\eh{i}|} \, , \qquad i = 1, 2, \dots, n-1 \, ,
  \end{aligned}
\end{equation}
and, finally, the $n$-th vector is defined as:
\begin{equation}
  \eh{n} = \star ( \fv{1} \wedge \star ( \fv{2} \wedge
  \star ( {\dots} \wedge \fv{n-1} ) \dots )) \, \Rightarrow \,  \fv{n} = \frac{\eh{n}}{|\eh{n}|} \, ,
\end{equation}
where $\wedge$ is the wedge product that can be thought as a
generalization of the cross product for vectors with dimensions
$n > 3$ \cite{Stoker}.  Observe that the wedge product of two vectors
is a bivector, hence the need for the Hodge star operator.

The generalized curvatures are given by:
\begin{equation}
  \chi_i = \dot {\mathbf{f}}_{i} \cdot \fv{i+1} \, ,
\end{equation}
with generalized frequencies:
\begin{equation}
  \omega_{\chi, i} = s' \chi_i = \fv{i}' \cdot \fv{i+1} \, .
\end{equation}
Defining $\F = [\fv{1}, \fv{2}, \dots, \fv{n}]$, the generalized
Frenet-Serret formulas become:
\begin{equation}
  \F' = \W{\chi} \, \F \, ,
\end{equation}
where
\begin{equation}
  \W{\chi} = 
  \begin{bmatrix}
    0 & \omega_{\chi, 1} & & 0 \\
    -\omega_{\chi, 1} & \ddots & \ddots & \\
    & \ddots & 0 & \omega_{\chi, n-1} \\
    0 & & -\omega_{\chi,n-1} & 0
  \end{bmatrix} .
\end{equation}
Finally, we note that Cartan's moving frames and attitude matrices
also immediately extend to $n$ dimensions.  In particular, observe
that, if $\bfb A$ is an attitude matrix with dimension $n$,
\eqref{eq:Aw0} returns a skew-symmetric matrix with dimension $n$.

\section{Geometrical Interpretation of the Voltage} 
\label{sec:geometry}

In \cite{freqgeom}, the author provides a geometrical interpretation
of electrical quantities in multi-phase circuits.  Limiting for
simplicity but without lack of generality the discussion to
three-phase circuits, the key assumption of \cite{freqgeom} is that
the voltage $\bfg v$ is a vector representing the time derivative of a
space curve.  According to the Faraday's law, this space curve has the
physical meaning of a magnetic flux vector, say
$\flux$:\footnote{Similarly, one can assume that the current in a
  three-phase line is the time derivative of a space curve, which by
  definition of electric current, has the meaning of an electric
  charge.}
\begin{equation}
  \label{eq:phi}
  \flux = - \bfg x \, ,  
\end{equation}
which, from the Faraday's law, leads to:
\begin{equation}
  \label{eq:v}
  \bfg v = -\bfp \flux = \bfp x \, .
\end{equation}
In \cite{freqfrenet}, the author utilizes \eqref{eq:v} to rewrite the
equations of the Frenet frame in terms of the voltage at a node of a
three-phase circuit, as follows.

From \eqref{eq:s} and \eqref{eq:v}, the magnitude of the voltage
vector is equivalent to the length $s$ of the trajectory of $\flux$.
Hence:
\begin{equation}
  \label{eq:sdot}
  s' = \bfm v \, ,
\end{equation}
which leads to conclude that $\bfm v$ is an invariant, as to be
expected.  Then, the following identities hold \cite{freqfrenet}:
\begin{equation}
  \label{eq:elecfrenet}
  \begin{aligned}
    &\T = \frac{\bfg v}{\bfm v}  \, , \qquad
    &\N = \frac{\sw \times \bfg v}{\bfm \omega \, \bfm v} \, , \qquad
    &\B = \frac{\sw}{\bfm \omega}  \, ,
  \end{aligned}
\end{equation}
where $\sw$ is the binormal vector before normalization:
\begin{equation}
  \label{eq:omega}
  \sw = \frac{\bfg v \times \bfp v}{|\bfg v|^2} \, .
\end{equation}

The curvature $\kappa$ and torsion $\tau$ that appear in the
Frenet-Serret equations \eqref{eq:frenet3} can be also expressed in
terms of the voltage, as follows:
\begin{equation}
  \label{eq:kappa2}
  \kappa = \frac{|\bfg v \times \bfp v|}{|\bfg v|^3} \, ,
\end{equation}
and:
\begin{equation}
  \label{eq:tau2}
  \tau = \frac{\bfg v \cdot \bfp v \times \bfpp v}{ |\bfg v \times \bfp v|^2 \, |\bfg v|^4} \, ,
\end{equation}
respectively.  From the latter two expressions and the identity
\eqref{eq:sdot}, one can also deduce that the azimuthal and torsional
frequencies that appear in \eqref{eq:Fw} are \cite{freqfrenet}:
\begin{equation}
  \label{eq:wkwt}
  \begin{aligned}
    \wk &= \bfm v \, \kappa = \bfm \omega \, , \\ 
    \wt &= \bfm v \, \tau \, .
  \end{aligned}
\end{equation}

It is important to note that differential geometry in general and the
Frenet frame in particular assume smooth curves.  This is a reasonable
assumption if the curve represents the trajectory of a point-mass in a
three-dimensional space (or even four-dimensional if one considers
space-time).  On the other hand, the fact that the voltage is a smooth
function of time does not necessarily always hold.  It can be argued
that instantaneous variations of the voltage are an approximation as,
in reality, parasite capacitive (and inductive) effects will always
make voltages (and currents) smooth state variables.  However, in
practice, very fast variations of the voltage complicate the
calculation of its time derivatives.  The case study discussed in
Section \ref{sub:case} shows that these discontinuities do not affect
the evaluation of the vectors of the Frenet frame and of the
quantities $\wk$ and $\wt$.

\subsection{Geometrical Interpretation of the Park Transform}
\label{sec:geopark}

We are now ready to present the main result of this work.  Let us
consider a balanced three-phase voltage:
\begin{equation}
  \label{eq:v1}
  \bfg v = v [\sin(\theta_a) \e{a} + \sin(\theta_a - \alpha) \e{b} +   
  \sin(\theta_a + \alpha) \e{c}] \, ,
\end{equation}
where $v$ and $\ta' = \wa$ are time-varying quantities.  Then, the
calculation of $(\T, \N, \B)$ based on \eqref{eq:elecfrenet} leads to:
\begin{equation}
  \label{eq:park1}
  \F = 
  \sqrt{\frac{2}{3}}
  \begin{bmatrix}
    \sin(\ta) & \sin(\ta - \alpha) & \sin(\ta + \alpha) \\
    \cos(\ta) & \cos(\ta - \alpha) & \cos(\ta + \alpha) \\
    \nicefrac{1}{\sqrt{2}} & \nicefrac{1}{\sqrt{2}} & \nicefrac{1}{\sqrt{2}}
  \end{bmatrix} \, ,
\end{equation}
which, comparing to \eqref{eq:park}, indicates that for a balanced
voltage and assuming $\tp = \theta_a$ the following identities hold:
\begin{equation}
  \label{eq:id1}
  \begin{aligned}
    \P(\wa) &= \F\, , &\quad \W{P} &= \W{F} \, , \\
    \wk &= \wa \, , &\quad \wt &= 0 \, .
  \end{aligned}
\end{equation}

Relevant special cases of \eqref{eq:v1} are:
\begin{itemize}
\item Stationary, balanced voltages with $v=\rm const.$ and
  $\wa = \wo$.  $\P(\wo)$ is the conventional Park transform utilized
  for transmission grid elements. 
\item For balanced synchronous machines, $v$ is the stator voltage and
  $\wa = \wr$ is the rotor angular speed.  $\P(\omega_r)$ represents
  the conventional synchronous machine Park transform $\P(\omega_r)$
  utilized in transient stability analysis.
\end{itemize}

The identities in \eqref{eq:id1} also indicate that, in general,
$\P(\wp) \ne \F$.  In particular, in the transient following a fault,
the frequency of the bus voltages varies from point to point of the
grid \cite{Divider}.  This implies that, in transient conditions, for
any bus voltage, $\P(\wp) = \F$ if and only if $\wp=\wk$.

The main conceptual difference between the Park transform and the
Frenet frame is that, for the former, $\wp$ has to be defined
\textit{a priori} and is, in general, a quantity detached from the
actual behavior of the quantities to which the transform is applied.
On the other hand, the Frenet frame is defined based on the
instantaneous values of the quantity to which it is applied and $\wk$
is obtained as a byproduct of the calculation of the frame itself.

Applying the Frenet frame to a three-phase voltage has the following
effect:
\begin{equation}
  \label{eq:vtnb2}
  \vtnb{v} = \F \, \vabc{v} =
  \begin{bmatrix}
    |\vabc{v}| \\ 0 \\ 0
  \end{bmatrix} .
\end{equation}
In fact:
\begin{equation}
  \T \cdot \vabc{v}
  = \frac{\vabc{v}}{|\vabc{v}|} \cdot \vabc{v}
  = |\vabc{v}| \, , 
\end{equation}
and, by construction, $\N \cdot \vabc{v} = \bfg 0$ and
$\B \cdot \vabc{v} = \bfg 0$.

Equation \eqref{eq:vtnb2} indicates that the Frenet frame
\textit{always} -- i.e., not only in balanced stationary conditions --
makes the original vector of voltage equivalent to a dc voltage.  The
``price'' of this transformation is that the time derivative of such a
voltage includes two terms, the conventional translation $\vtnb{v}'$
and the rotation $\W{F} \vtnb{v}$.  In the other way round, thus, one
can view a dc voltage as a quantity referred to a Frenet frame for
which $\W{F} = \bfg 0$, which is, in effect, the condition satisfied
by straight lines.

In conclusion, at every instant, a three-phase voltage $\bfg v$ is
fully characterized by three scalar quantities, namely $\bfm v$, $\wk$
and $\wt$.  Again, this result does not apply only to balanced
stationary conditions, but always hold.  Similarly, $n$-dimensional
voltage vectors are fully characterized by $\bfm v$ and $n-1$
generalized frequencies, namely
$\omega_{\chi, 1}, \dots, \omega_{\chi, n-1}$.

\subsection{Beyond Balanced Conditions}

The previous section shows that, for balanced stationary conditions
and if $\wp = \omega_a$, the Park transform and the Frenet frame
coincide.  In these conditions, in fact, the two transforms differ at
most by a phase shift.  This occurs if the Park transform reference
angle is $\tp = \theta_a + \theta_{{\scriptscriptstyle \rm P}, 0}$.
The differences -- and generality -- of the Frenet frame with respect
to the Park transform, however, manifests in unbalanced and/or
transient conditions.

The first difference is that the Frenet frame always \textit{follows}
the curve and defines, at each instant, a set of orthogonal
coordinates that have a specific meaning for the trajectoy itself
(namely, tangent, normal and binormal vectors).  For this reason,
$\wk$ (and $\wt$) is a \textit{byproduct} of the Frenet frame, not an
\textit{arbitrary choice} as it is $\wp$ for the Park transform.  The
second difference is that the Frenet frame does not require setting an
external and, again, arbitrary, reference phase angle.

The differences above have relevant consequences.  For example, all
trajectories that lay in the plane $(\T, \N)$ have $\wt = 0$.
Reference \cite{freqfrenet} shows that voltages that are stationary
unbalanced, balanced with harmonics, and balanced in transient
conditions show $\wt = 0$, whereas stationary voltages with unbalanced
harmonic content show $\wt \ne 0$.  On the other hand, since the Park
transform does not follow the trajectory described by the voltage, it
shows a nonnull $o$-axis component in unbalanced conditions.  More
importantly, since the frequency $\wp$ is not necessarily related to
the time evolution of the voltage, if $\wp \ne \omega_a$, Park
transformed $dq$-axis components can vary in time also in stationary
conditions.  This situation is common in power systems, where, due to
the droop and deadbands of primary frequency controllers, the
frequency of the voltages can be different (even if just slightly)
from the synchronous reference $\omega_o$.  Instead, the Frenet frame
satisfies by construction the condition $\wk = \omega_a$ in stationary
balanced conditions.

In turn, the Park transform has a relevant geometrical meaning only in
the balanced stationary case.  In all other conditions, the components
of the Park transform are simply projections on an arbitrary
(\textit{arbitrary} in the sense that the coordinates are not related
in any way to the trajectory described by the voltage) sets of
time-varying coordinates.  The geometrical meaning of the Frenet frame
and the and properties of the Frenet transformed quantities and
invariants, on the other hand, are always the same.  Based on this
observation, the Frenet transform can be seen as a generalization of
the Park transform.

Sections \ref{sub:unbalanced} to \ref{sub:case} further illustrate,
through numerical examples, all points discussed above.

\subsection{Interconnection of Voltage Frenet Frames to the Grid}

The ``locality'' of the Frenet frame of each individual device has to
be conciliated with the rest of the grid.  To be able to study the
interaction of the devices and, ultimately, to study and simulate the
dynamics of the system, it is thus necessary to have a mechanism to
convert the local set of coordinates $(\T, \N, \B)$ to the ``system''
coordinates.

In Section \ref{sec:park}, we have described the common way with which
synchronous machines are interfaces to the grid, namely using a
rotation in the plane $dq$ through a cylindrical coordinate change as
in \eqref{eq:C}.  Since the Frenet frame provides a systematic way to
define the ``local'' frequency of a device as its azimuthal frequency,
a change of coordinate can be done, in effect, with any combination of
transforms.  For example, assuming that the network is referred to the
Park transform $\P(\ws)$, the change of coordinate from a local Frenet
frame of a device to the network is given by:
\begin{equation}
  \vdqo{v}^N =
  \P(\ws) \, \Tr{F} \, \vtnb{v}^D =
  \bfb \Psi(\ws) \, \vtnb{v}^D \, ,
\end{equation}
where $N$ denotes a quantity in the network reference frame as in
\eqref{eq:G2N}, and $D$ denotes a quantity expressed using a local
reference frame.  Then, the rotation of the attitude matrix
$\bfb \Psi$ is given by:
\begin{equation}
  \label{eq:PwFw}
  \begin{aligned}
    \W{\Psi}
    &= \big (\P(\ws) \, \Tr{F} \big )' \, \big (\widehat{ \P(\ws) \, \Tr{F} }\big ) \\
    &= \big (\P'(\ws) \, \Tr{F} + \P(\ws) \, {\Tr{F}}' \big ) \, \F \, \Tr{P}(\ws) \\
    &= \P'(\ws) \, \Tr{F} \, \F \, \Tr{P}(\ws) +
    \P(\ws) \, \Tr{F}' \, \F \, \Tr{P}(\ws) \, .
  \end{aligned}
\end{equation}
Finally, recalling that $\Tr{F} \, \F = \bfb I$, $\Wt{F} = \F \, \Tr{F}'$,
and $\Wt{F} = -\W{F}$, \eqref{eq:PwFw} becomes:
\begin{equation}
  \label{eq:Wpsi}
  \W{\Psi} = \W{P} - \bfb \Psi(\ws) \, \W{F} \, \Tr{\Psi}(\ws) \, ,
\end{equation}
which is a skew-symmetric matrix in the form of
\eqref{eq:Aw1}.\footnote{Note that \eqref{eq:Wpsi} holds in general
  for any product of attitude matrices.  In fact, if $\A$ and $\bfb B$
  are attitude matrices of same order, then:
  \begin{equation*}
    \bfb E = \A \, \Tr{B}
    \quad \Rightarrow \quad
    \W{E} = \W{A} - \bfb E \, \W{B} \Tr{E} \, .
  \end{equation*}}
  %
  %
%
Matrix $\bfb \Psi$ is the Frenet frame equivalent to the cylindrical
frame $\C$ defined in \eqref{eq:C}.  Or, equivalently, $\bfb \Psi$ is
the generalization of $\C$.  In the same vein, the expression
\eqref{eq:Cw} represents the rotation matrix $\W{C}$, of which
$\W{\Psi}$ is the generalization.

\section{Case Study}
\label{sec:examples}

This section illustrates the theoretical results above through a
variety of examples, including a simulation based on the IEEE 39-bus
system.  In Sections \ref{sub:balanced} to \ref{sub:case}, the
voltages $\vabc{v}$ are assumed to be three-dimensional vectors on the
basis:
\begin{equation}
  \label{eq:basis0}
  \begin{aligned}
    \ei{a} &= (1, 0, 0) \, , \quad 
    \ei{b} &= (0, 1, 0) \, , \quad
    \ei{c} &= (0, 0, 1) \, . 
  \end{aligned}
\end{equation}
The last example (Section \ref{sub:6phase}) considers a multi-phase
system with $n=6$.  In all cases, $\wo = 2\pi \, 60$ rad/s and
$\P = \P(\wo)$.  In the figures below, the trajectories of the voltage
are in pu with respect to a base of 15 kV and $\wk$ and $\wt$ are in
pu with respect to $\wo$.

\begin{figure*}[ht!]
  \centering
  \resizebox{0.97\linewidth}{!}{\includegraphics{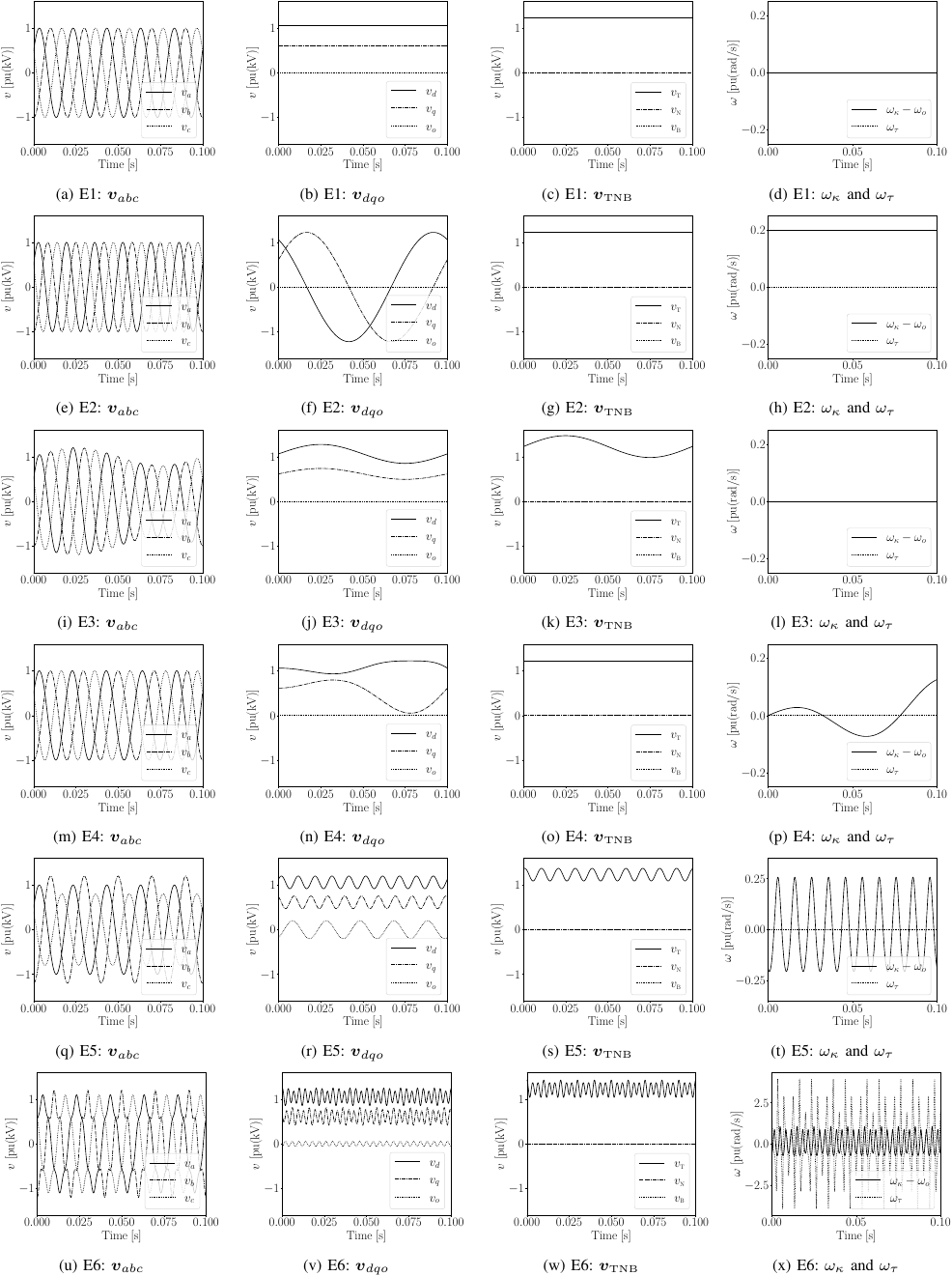}}
  \caption{Examples E1-E6.}
  \label{fig:ex}
\end{figure*}

\subsection{Balanced Three-Phase AC Voltages}
\label{sub:balanced}

The examples presented in this section utilize a balanced three-phase
AC voltages in the form:
\begin{equation}
  \begin{aligned}
    \label{eq:balanced}
    v_a &= V \sin (\w{a} \, t  + \th{a,0}) \, , \\ 
    v_b &= V \sin (\w{a} \, t + \th{a,0} + \alpha) \, , \\ 
    v_c &= V \sin (\w{a} \, t + \th{a,0} - \alpha) \, ,
  \end{aligned}
\end{equation} 
where $\th{a,0} = \pi/6$ rad.  The following four cases are
considered:
\begin{itemize}
\item E1: $\wp = \wa = \wo$, $V = 15$ kV. 
\item E2: $\wp = \wo$, $\wa = 1.2 \, \wo$ rad/s, $V = 15$ kV. 
\item E3: $\wp = \wa = \wo$, $V = 15 + 3 \sin(0.2\wo t)$ kV. 
\item E4: $\wp = \wo$, $\w{a} = \wo + 2 \pi \sin(20 \pi t)$ rad/s,
  $V = 15$ kV.
\end{itemize}

The voltage $\vabc{v}$, $\vdqo{v}$, $\vtnb{v}$ as well as the angular
frequencies $\wk$ and $\wt$ for the examples E1-E4 are shown in
Figs.~\ref{fig:ex}.(a)-(p).

Example E1 is the standard balanced, stationary case, for which the
Park transform and the Frenet frame substantially coincide, except for
the fact that the Frenet frame is independent from the angular
position $\th{a,0}$.  This is a consequence of the fact that the
Frenet frame follows the voltage locally, whereas the Park transform
follows an independent reference.

Example E2 shows one of the critical issues (especially, in the
context of state estimation) of the Park transform: the fact that if
$\wp \ne \w{a}$ then, the Park $dq$-axis components oscillates at
frequency $\wa -\wp$, even if the amplitude of the measured signal is
perfectly stationary \cite{Paolone:2020}.  As expected, on the other
hand, the Frenet frame returns a steady-state voltage magnitude $v$
and azimuthal frequency $\wk$.

Examples E3 and E4 represent dual scenarios.  In E3, the voltage
magnitude is time varying and the frequency is constant whereas, in
E4, the voltage magnitude is constant and the angular frequency
$\w{a}$ is time varying.  The Park transform is unable to distinguish
between these two situations, once again because of the constant
external reference angular speed.  On the other hand, the Frenet frame
separates the effects of the variations of the voltage magnitude and
of the angular frequency.  It is relevant to observe that this
property holds also for signals with time-varying voltage magnitude
and angular frequency, which is the typical scenario in the first
seconds of power system transients following a large disturbance (see
the case study in Section \ref{sub:case}).

\subsection{Unbalanced Three-Phase AC Voltages}
\label{sub:unbalanced}

This section discusses the effect of unbalanced conditions.  Let us
define the following unbalanced three-phase AC voltages (example E5):
\begin{equation}
  \begin{aligned}
    \label{eq:unbalanced}
    v_a &= V \sin (\wo \, t  + \th{a,0}) \, , \\ 
    v_b &= 1.2 V \sin (\wo \, t + \th{a,0} + \alpha) \, , \\ 
    v_c &= 0.8 V \sin (\wo \, t + \th{a,0} - \alpha) \, ,
  \end{aligned}
\end{equation} 
where $\th{a,0} = \pi/6$ rad.  The voltage $\vabc{v}$, $\vdqo{v}$,
$\vtnb{v}$ as well as the angular frequencies $\wk$ and $\wt$ for the
voltage vector \eqref{eq:unbalanced} are shown in
Figs.~\ref{fig:ex}.(q)-(t).

The unbalanced conditions give birth to a non-null zero-sequence
component in the vector $\vdqo{v}$.  However, the curve associated
with these conditions is still a plane curve, as confirmed by the fact
that $\wt = 0$.  In turn, the trajectory described by the voltage is
an ellipse rather than a circle.  Then, to describe the curve it
suffices to know only two quantities, not three as suggested by the
Park transforms.  This situation is properly captured by the Frenet
frame: both the magnitude of the voltage and the azimuthal frequency
are periodic, which reflects the fact that in ellipse, both the radius
and the curvature are not constant and repeat periodically at every
full turn.

\subsection{Three-Phase AC Voltages with Harmonics}
\label{sub:harmonics}

This example (E6) describes the effect of unbalanced harmonics on
voltages transformed with the Park transform and the Frenet frame.  We
consider the following voltage vector:
\begin{equation}
  \begin{aligned}
    \label{eq:harmonics}
    v_a &= V \, \{\sin (\wo \, t  + \th{a,0}) + 0.1 \,
    \sin [5(\wo \, t  + \th{a,0})]\} \, , \\ 
    v_b &= V \, \{\sin (\wo \, t  + \th{b,0}) + 0.2 \,
    \sin [5(\wo \, t  + \th{b,0})]\} \, , \\ 
    v_c &= V \, \{\sin (\wo \, t  + \th{c,0}) + 0.1 \,
    \sin [5(\wo \, t  + \th{c,0})]\} \, ,
  \end{aligned}
\end{equation}
where $V = 15$ kV, $\th{a,0} = \pi/6$ rad, $\th{b,0} = \pi/6 - \alpha$
rad, and $\th{c,0} = \pi/6 + \alpha$ rad.  The voltage $\vabc{v}$,
$\vdqo{v}$, $\vtnb{v}$ as well as the angular frequencies $\wk$ and
$\wt$ for the voltage vector \eqref{eq:harmonics} are shown in
Figs.~\ref{fig:ex}.(u)-(x).

Unbalanced harmonics leads to a time-varying torsional frequency $\wt$
(see also the examples included in \cite{freqfrenet}).  It is relevant
to note that, from the Park transform point of view, both E5 and E6
lead to a periodic $v_o$.  Of course, the frequency of the
oscillations allows distinguishing between E5 (unbalance voltages at
the fundamental frequency) and E6 (harmonic content).  However, only
the Frenet frame allows interpreting correctly the ``non-planar''
nature of E6.

\subsection{Power System Transient}
\label{sub:case}

\begin{figure}[ht!]
  \centering
  \resizebox{0.98\linewidth}{!}{\includegraphics{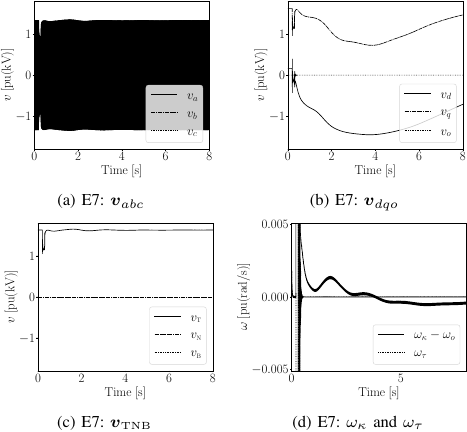}}
  \caption{Example E7: Transient behavior of the voltage and frequency
    at bus 26 of the IEEE 39-bus system following a fault at bus 4.}
  \label{fig:sims}
  \vspace{-1mm}
\end{figure}

This last example considers the IEEE 39-bus system.  The model
utilized in the simulation below is provided as an application example
with DIgSILENT PowerFactory.  This consists of a detailed EMT
three-phase model.  Beside machine and control dynamics of the
original IEEE 39-bus benchmark system, the model includes the
electromagnetic dynamics of transmission lines, transformers and
synchronous machines.  The contingency is a three-phase fault at bus
4.  The fault occurs at $t= 0.2$ s and is cleared at $t = 0.3$ s.  The
time step of the numerical integration is $0.01$ ms.  Figure
\ref{fig:sims} shows the trajectories of $\vabc{v}$, $\vdqo{v}$ and
$\vtnb{v}$ at bus 24 in pu with respect to a base of $220$ kV and of
the azimuthal and torsional frequencies in pu with respect to
$\wo = 2\pi \, 60$ rad/s.

The untransformed voltage vector $\vabc{v}$ carries ``too much''
information for the time scale of the simulation (several seconds),
which is the typical one utilized to study electromechanical
transients after a fault.  This fully justifies the common practice to
utilize a RMS model for the study of these transients.  The
Park-transformed voltage $\vdqo{v}$ suffers of the decoupling between
the reference angular frequency utilized in $\P(\wo)$.  The slip
between the actual frequency at the bus and $\wo$ makes the
interpretation of the behavior of the $dq$-axis components $v_d$ and
$v_q$ not intuitive.  The Frenet frame, on the other hand, provides
straightforward information.  The tangent component
$v_{\scriptscriptstyle \rm T}$ is the time envelope of $\vabc{v}$.
Then $\wk$ shows the behavior of the local frequency at the bus.
Finally, $\wt$, which is null except during the fault, indicates that
the system does not include unbalanced harmonics.

On a practical note, the accuracy of the evaluation of the vectors of
the Frenet frame as well as of the geometric invariants depends
exclusively to the precision with which the time derivatives of the
voltage can be evaluated.  In fact, $\T$ is the normalized three-phase
ac voltage $\bfg v_{abc}$, whereas $\N$ and $\B$ are obtained based on
the first time derivative of $\bfg v_{abc}$, see \eqref{eq:elecfrenet}
and \eqref{eq:omega}.  Similarly $\omega_{\kappa}$ requires only the
calculation of first time derivatives, whereas $\omega_{\tau}$
requires also second time derivatives, see \eqref{eq:kappa2} and
\eqref{eq:tau2}.  In this case study, the derivatives were obtained by
numerically differentiating the sampled phase voltages at bus 26 and
then removing noise with discrete butterworth and low-pass filters.

Finally, it is relevant to remark that, for a balanced case as the one
discussed in this section,
$v_{\scriptscriptstyle \rm N} = \sqrt{v_d^2 + v_q^2}$ and
$\wk = \wo + \frac{d}{dt} \arctan(v_q, v_d)$.  Moreover, the drift of
the $dq$-axis components can be compensated by using the frequency of
the center of inertia rather than $\wo$ in matrix $\P$.  However, the
center of inertia is a quantity that can be calculated only in a
computer simulation, not in practice.  And, more importantly and as
shown from examples E3 to E6, the Park transform does not
\textit{always} provide an easy way to interpret the results.  The
feature of decomposing the original voltage vector into meaningful
quantities (invariants) is, by construction, specific only of the
Frenet frame and constitutes thus its main advantage with respect to
any other transforms.

\vspace{-2mm}
\subsection{Balanced Six-Phase Voltage}
\label{sub:6phase}

The Gram-Schmidt process described in Section \ref{sub:ndim} is
illustrated using a balanced six-phase system.  Consider the following
voltage vector:
\begin{equation}
  \label{eq:v6}
  \bfg v = V \sum_{h=1}^6 \sin \big (\wo t - (h-1) \beta \big ) \e{h} \, ,
\end{equation}
The first two orthonormal vectors of the Frenet frame are:
\begin{equation}
  \begin{aligned}
    \eh{1} &= \frac{\bfg v}{\bfm v} \, \, \Rightarrow \, \,
    \fv{1} = \frac{1}{\sqrt{3}}
    \sum_{h=1}^6 \sin \big (\wo t - (h-1) \beta \big ) \e{h} \, , \\
    \eh{2} &= \frac{\bfg v'}{\bfm{v}^2} \, \, \Rightarrow \, \,
    \fv{2} = \frac{1}{\sqrt{3}}
    \sum_{h=1}^6 \cos \big (\wo t - (h-1) \beta \big ) \e{h} \, .
  \end{aligned}
\end{equation}
Then, note that $\bfg v'' = - \wo^2 \bfg v$ and
$\bfg v''' = -\wo^2 \bfg v'$, etc.  This indicates that the balanced
six-phase voltage \eqref{eq:v6} is in effect a plane curve. The only
nonnull remaining vector of the basis is thus the normalized vector
perpendicular to both $\fv{1}$ and $\fv{2}$, namely:
\begin{equation}
  \eh{3}
  = \star (\fv{1} \wedge \fv{2}) \, \,  \Rightarrow \, \,
  \fv{3} = \frac{1}{\sqrt{6}} \sum_{h=1}^6 \ei{i} \, .
\end{equation}

The obtained Frenet frame $\F = [\fv{1}, \fv{2}, \fv{3}]$ is exactly
the same 6-to-$dqo$ Park transform matrix proposed in
\cite{GABER1988203}.  The procedure utilized in \cite{GABER1988203},
however, is more involved than the one proposed here as it utilizes
the definition of groups.  We have thus obtained that for a stationary
balanced $n$-dimensional voltage vector the generalized Park transform
is equivalent to the generalized Frenet frame.  However, as expected,
the Frenet frame is more general as it does not require the voltage to
be balanced or stationary.  Finally, we observe that
$\omega_{\chi, 1} = \wo$ and $\omega_{\chi, i} = 0$ for
$i = 2, \dots, 6$, which confirms that the curve lays on a plane.

\vspace{-1mm}
\section{Conclusions}
\label{sec:conc}

The paper presents a geometrical interpretation of the Park transform,
its time derivative and its generalization based on the Frenet frame
and Cartan's moving frames.  The Frenet frame appears particularly
relevant for the transient stability analysis of power systems for
various reasons, as follows.

\begin{itemize}
\item The Frenet frame returns a set of geometric invariants which are
  as many as the phases of the circuit.  These quantities represent
  fully and unequivocally the transient conditions of the voltage (or
  current) under consideration.
\item The Frenet frame provides a natural phase angle reference as
  well as intrinsic angular frequency (namely the azimuthal frequency)
  for the voltage without the need for an external reference or a
  device that links to an external reference, such as the phase-locked
  loops.
\item Cartan's moving frame approach provides a systematic and general
  way to link the local Frenet frames to a common reference.  The
  equations of such interfaces are duly provided in this work.
\item The paper also shows how the Frenet frame can be extended to any
  number of phases and provides the steps required to calculate the
  generalized coordinates and curvatures for an arbitrary multi-phase
  circuit.
\end{itemize}

In the case study section, a variety of examples support the theory
and show how the proposed approach solves the many idiosyncrasies of
the Park transform.

The ability to identify the angular frequency of devices that do not
have a rotor appears particularly promising for the study of
converter-interfaced generation.  Another relevant aspect is the
practicality of the implementation of the proposed technique for
on-line applications, such as control and dynamic state-estimation.
These topics will be the focus of future work.


\begin{thebibliography}{10}
\providecommand{\url}[1]{#1}
\csname url@samestyle\endcsname
\providecommand{\newblock}{\relax}
\providecommand{\bibinfo}[2]{#2}
\providecommand{\BIBentrySTDinterwordspacing}{\spaceskip=0pt\relax}
\providecommand{\BIBentryALTinterwordstretchfactor}{4}
\providecommand{\BIBentryALTinterwordspacing}{\spaceskip=\fontdimen2\font plus
\BIBentryALTinterwordstretchfactor\fontdimen3\font minus
  \fontdimen4\font\relax}
\providecommand{\BIBforeignlanguage}[2]{{%
\expandafter\ifx\csname l@#1\endcsname\relax
\typeout{** WARNING: IEEEtran.bst: No hyphenation pattern has been}%
\typeout{** loaded for the language `#1'. Using the pattern for}%
\typeout{** the default language instead.}%
\else
\language=\csname l@#1\endcsname
\fi
#2}}
\providecommand{\BIBdecl}{\relax}
\BIBdecl

\bibitem{5055275}
R.~H. Park, ``Two-reaction theory of synchronous machines generalized method of
  analysis -- {Part I},'' \emph{Transactions of the American Institute of
  Electrical Engineers}, vol.~48, no.~3, pp. 716--727, 1929.

\bibitem{FDF:2020}
F.~{Milano} and {\'A}.~{Ortega}, \emph{Frequency Variations in Power Systems:
  Modeling, State Estimation, and Control}.\hskip 1em plus 0.5em minus
  0.4em\relax Hoboken, NJ: Wiley, 2020.

\bibitem{Iravani:2010}
A.~Yazdani and R.~Iravani, \emph{Voltage-Sourced Converters in Power Systems:
  Modeling, Control, and Applications}.\hskip 1em plus 0.5em minus 0.4em\relax
  Hoboken, NJ: Wiley, 2010.

\bibitem{GABER1988203}
A.~Z. Gaber and L.~P. Singh, ``Investigation of symmetries inherent in
  synchronous machines and development of a generalized {Park}'s transformation
  based solely upon symmetries,'' \emph{Electric Power Systems Research},
  vol.~15, no.~3, pp. 203--213, 1988.

\bibitem{4454446}
E.~Levi, ``Multiphase electric machines for variable-speed applications,''
  \emph{IEEE Transactions on Industrial Electronics}, vol.~55, no.~5, pp.
  1893--1909, 2008.

\bibitem{Steinmetz:1897}
C.~P. Steinmetz and E.~J. Berg, \emph{Theory and Calculation of Alternating
  Current Phenomena}.\hskip 1em plus 0.5em minus 0.4em\relax New York, US:
  W.~J.~Johnston Co., 1897.

\bibitem{Das:2015}
J.~C. Das, \emph{Power System Harmonics and Passive Filter Designs}.\hskip 1em
  plus 0.5em minus 0.4em\relax Hoboken, NJ: IEEE Press -- John Wiley \& Sons,
  2015.

\bibitem{Fortescue:1918}
C.~L. Fortescue, ``Method of symmetrical co-ordinates applied to the solution
  of polyphase networks,'' \emph{Transactions of the American Institute of
  Electrical Engineers}, vol.~37, no.~2, pp. 1027--1140, June 1918.

\bibitem{Ku:1929}
Y.~H. Ku, ``Transient analysis of {A-C.}~machinery,'' \emph{Transactions of the
  American Institute of Electrical Engineers}, vol.~48, no.~3, pp. 707--714,
  July 1929.

\bibitem{Clarke:1943}
E.~Clarke, \emph{Circuit Analysis of AC Power Systems -- Volume I: Symmetrical
  and Related Components}, ser. General Electric Series.\hskip 1em plus 0.5em
  minus 0.4em\relax New York, US: J. Wiley \& Sons, 1943.

\bibitem{871734}
A.~M. {Stankovi{\'c}} and T.~{Aydin}, ``Analysis of asymmetrical faults in
  power systems using dynamic phasors,'' \emph{IEEE Transactions on Power
  Systems}, vol.~15, no.~3, pp. 1062--1068, Aug. 2000.

\bibitem{6687278}
C.~{Liu}, A.~{Bose}, and P.~{Tian}, ``Modeling and analysis of {HVDC} converter
  by three-phase dynamic phasor,'' \emph{IEEE Transactions on Power Delivery},
  vol.~29, no.~1, pp. 3--12, Feb. 2014.

\bibitem{Hancock:1974}
N.~N. Hancock, \emph{Matrix Analysis of Electrical Machinery}, 2nd~ed., ser.
  General Electric Series.\hskip 1em plus 0.5em minus 0.4em\relax New York, US:
  Pergamon Press, 1974.

\bibitem{Stoker}
J.~J. Stoker, \emph{Differential Geometry}.\hskip 1em plus 0.5em minus
  0.4em\relax New York: Wiley-Interscience, 1969.

\bibitem{ONeill}
B.~{O'Neill}, \emph{Elementary Differential Geometry}.\hskip 1em plus 0.5em
  minus 0.4em\relax London: Academic Press, 1966.

\bibitem{2003lapierre}
L.~Lapierre, D.~Soetanto, and A.~Pascoal, ``Nonlinear path following with
  applications to the control of autonomous underwater vehicles,'' in
  \emph{42nd IEEE International Conference on Decision and Control}, vol.~2,
  2003, pp. 1256--1261.

\bibitem{2010werling}
M.~Werling, J.~Ziegler, S.~Kammel, and S.~Thrun, ``Optimal trajectory
  generation for dynamic street scenarios in a {Frenet} frame,'' in \emph{IEEE
  International Conference on Robotics and Automation}, 2010, pp. 987--993.

\bibitem{freqfrenet}
F.~Milano, G.~Tzounas, I.~Dassios, and T.~K{\"e}r{\c{c}}i, ``Applications of
  the {Frenet} frame to electric circuits,'' \emph{IEEE Transactions on
  Circuits and Systems I: Regular Papers}, vol.~69, no.~4, pp. 1668--1680,
  2022.

\bibitem{freqgeom}
F.~Milano, ``A geometrical interpretation of frequency,'' \emph{IEEE
  Transactions on Power Systems}, vol.~37, no.~1, pp. 816--819, 2022.

\bibitem{4450060503}
J.~L. Willems, ``Mathematical foundations of the instantaneous power concepts:
  A geometrical approach,'' \emph{European Transactions on Electrical Power},
  vol.~6, no.~5, pp. 299--304, 1996.

\bibitem{1308315}
X.~Dai, G.~Liu, and R.~Gretsch, ``Generalized theory of instantaneous reactive
  quantity for multiphase power system,'' \emph{IEEE Transactions on Power
  Delivery}, vol.~19, no.~3, pp. 965--972, 2004.

\bibitem{5316097}
H.~{Lev-Ari} and A.~M. {Stankovi{\'c}}, ``Instantaneous power quantities in
  polyphase systems -- {A} geometric algebra approach,'' in \emph{IEEE Energy
  Conversion Congress and Exposition}, 2009, pp. 592--596.

\bibitem{IPT}
H.~Akagi, E.~H. Watanabe, and M.~Aredes, \emph{Instantaneous Power Theory and
  Applications to Power Conditioning}, 2nd~ed.\hskip 1em plus 0.5em minus
  0.4em\relax New York: Wiley IEEE Press, 2017.

\bibitem{math9111295}
F.~G. Montoya, R.~Ba{\~n}os, A.~Alcayde, F.~M. Arrabal-Campos, and J.~Rold{\'
  a}n-P{\'e}rez, ``Vector geometric algebra in power systems: An updated
  formulation of apparent power under non-sinusoidal conditions,''
  \emph{Mathematics}, vol.~9, no.~11, 2021.

\bibitem{2016barry}
N.~Barry, ``The application of quaternions in electrical circuits,'' in
  \emph{2016 27th Irish Signals and Systems Conference (ISSC)}, 2016, pp. 1--9.

\bibitem{2018ishihara}
V.~d.~P. Brasil, A.~de~Leles Ferreira~Filho, and J.~Y. Ishihara, ``Electrical
  three phase circuit analysis using quaternions,'' in \emph{18th International
  Conference on Harmonics and Quality of Power (ICHQP)}, 2018, pp. 1--6.

\bibitem{2015talebi}
S.~P. Talebi and D.~P. Mandic, ``A quaternion frequency estimator for
  three-phase power systems,'' in \emph{IEEE International Conference on
  Acoustics, Speech and Signal Processing (ICASSP)}, 2015, pp. 3956--3960.

\bibitem{Fabozzi:2011}
D.~Fabozzi and T.~{Van Cutsem}, ``On angle references in long-term time-domain
  simulations,'' \emph{IEEE Transactions on Power Systems}, vol.~26, no.~1, pp.
  483--484, Feb. 2011.

\bibitem{Needham:2021}
T.~Needham, \emph{Visual Differential Geometry and Forms: A Mathematical Drama
  in Five Acts}.\hskip 1em plus 0.5em minus 0.4em\relax Princeton, NJ:
  Princeton University Press, 2021.

\bibitem{Divider}
F.~{Milano} and {\'A}.~{Ortega}, ``Frequency divider,'' \emph{IEEE Transactions
  on Power Systems}, vol.~32, no.~2, pp. 1493--1501, 2017.

\bibitem{Paolone:2020}
{M. Paolone {\em et al.}}, ``Fundamentals of power systems modelling in the
  presence of converter-interfaced generation,'' \emph{Electric Power Systems
  Research}, vol. 189, p. 106811, 2020.

\end{thebibliography}


\begin{IEEEbiography}
  [{\includegraphics[width=1in, height=1.25in, clip,
    keepaspectratio]{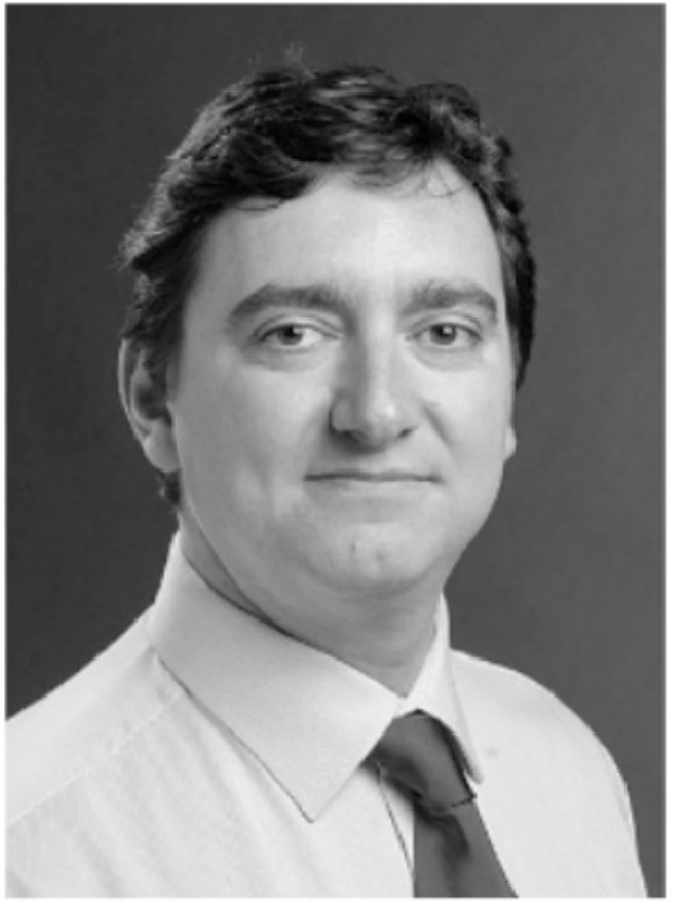}}] {Federico Milano} (F'16) received
  from the Univ. of Genoa, Italy, the ME and Ph.D.~in Electrical
  Engineering in 1999 and 2003, respectively.  From 2001 to 2002 he
  was with the University of Waterloo, Canada, as a Visiting Scholar.
  From 2003 to 2013, he was with the University of Castilla-La Mancha,
  Spain.  In 2013, he joined the University College Dublin, Ireland,
  where he is currently a full professor.  He is also Chair of the
  IEEE Power System Stability Controls Subcommittee, IET Fellow, IEEE
  PES Distinguished Lecturer, and Co-Editor in Chief of the IET
  Generation, Transmission \& Distribution.  His research interests
  include power system modeling, control and stability analysis.
\end{IEEEbiography}

\vfill

\end{document}